\documentclass[letter, 12pt]{article}
\usepackage[T1]{fontenc}
\usepackage{graphicx}
\usepackage[ansinew]{inputenc}
\usepackage{latexsym}
\usepackage{amsthm}
\usepackage{amssymb}
\usepackage{amsmath}
\usepackage{multirow}
\usepackage{verbatim}

\newtheorem{thm}{Theorem}[section]
\newtheorem{prop}[thm]{Proposition}
\newtheorem{Lem}{Lemma}[section]

\theoremstyle{remark}
\newtheorem{rem}{Remark}[section]
\theoremstyle{definition}
\newtheorem{defn}{Definition}[section]
\theoremstyle{plain}
\newenvironment{dem}{\textbf{Proof.\ }}{\hfill$\Box$}

\usepackage{geometry}
\title{\textbf{A class of copulae associated with Brownian motion processes and their maxima}}

\begin{document}

	\maketitle


{Michel Ad\`es$^{1*}$,  Matthieu Dufour$^1$, Serge B. Provost$^{2}$, Marie-Claude Vachon$^1$}\\

\noindent $^1${\textit{D\'epartement de Math\'ematiques, Universit\'e du Qu\'ebec \`a Montr\'eal, Qu\'ebec, Canada}}\\
\noindent $^2${\textit{Department of Statistical and Actuarial Sciences,  The University of Western Ontario, London, Canada}}\\
\noindent $*$ Corresponding author:  ades.michel@uqam.ca

\newpage

\begin{abstract}

\noindent The main objective of this paper consists in creating a new class of copulae from various joint distributions occurring in connection with certain Brownian motion processes. We focus our attention on the distributions of univariate Brownian motions having a drift parameter and their maxima and on correlated bivariate Brownian motions by considering the maximum value of one of them. 
The copulae generated therefrom
 and their associated density functions are explicitly given  as well as graphically represented.
\end{abstract}
\

\noindent\textbf{Keywor{\rm{d}}s}: Brownian motion, Copulas, Correlated Brownian processes, Dependence, Two-dimensional Brownian motion. 
\

\

\section{Introduction}
This section first presents useful background information on Brownian motion $(\cal{BM})$. Then, copulae are defined and relevant related results are provided.

In 1918, the mathematician Norbert Wiener gave a rigorous formulation of Brownian motion and established its existence, which explains why the alternative name, Wiener process, is also in use.
 $\cal{BM}$  is utilized in various fiel{\rm{d}}s of scientific investigation such as 
Economics, Biology, Communications Theory, Business Administration, and Quantitative Finance. 
For instance, as pointed out by Chuang (1994), distributional results for $\cal{BM}$ can also be utilized for pricing contingent claims with barriers on price processes; Cao (2017) made use of correlated Brownian motions to solve an optimal investment-reinsurance problem.

Let $\{W_t\}_{t\geq 0}$ and $W_T$ represent the standard $\cal{BM}$ process and its terminal value, $M_t=\max_{0\leq s\leq t}W_s$ and  $M_{(s,\,t)}=\max \{W_u,\,s\leq u\leq t\}$. We shall consider the joint distributions of 
\begin{enumerate}
\item $W_t$ and its maximum $M_t$\,,
\item $W_T$ and $M_t$\,,
\item $W_T$ and $M_{(s,\,t)}$\,,
\end{enumerate}
which have previously been studied by Harrison (1985), Chuang (1996) and Lee (2003), among others.

Some further related results are available in the statistical literature. For example, representations of the joint density function of a  $\cal{BM}$  process and its minimum and maximum, which are given for instance in Borodin and Salminen (2002), were shown to be convergent by Choi and Roh (2013). Upper and lower boun{\rm{d}}s for the distribution of the maximum of a two-parameter  $\cal{BM}$  process were obtained by Caba{\~n}a and Wschebor (1982). Vardar-Acara et al. (2013) provided explicit expressions for the correlation between the supremum and the infimum of a $\cal{BM}$ with drift.
Kou and Zhong (2016) studied the first-passage times of two-dimensional  $\cal{BM}$ processes. Haugh (2004) explained how to generate correlatated Brownian motions and  points out some applications involving security pricing and porfolio evaluation.

We now review  some basic definitions and theorems in connection with copulae. Additional results are available from several authors including Cherubini et al. (2004, 2012), Denuit et al. (2005), Joe (1997), Nelsen (2006), and Sklar (1959).

The main idea behind copulae is that the joint distribution of two or more random variables can be expressed in terms of their marginal distributions and a certain correlation structure. Copulae enable one to separate the effect due to the dependence between the variables from contribution of each of the marginal variables. We focus on the two-dimensional case in this paper. In this framework, a copula function is a bivariate distribution defined on the unit square $\mathbf{I}^2=[0,\,1]^2$ with uniformly distributed margins. Formally, we have:

\begin{defn}
A function $C:\mathbf{I}^2\mapsto \mathbf{I}$ is a bivariate copula if it satisfies the following properties:
\begin{enumerate}
\item For every $y,\,w\in \mathbf{I}$,
\begin{eqnarray*}
C(y,\,1)&=&y\ \textrm{\  and \ }\  C(1,\,w)=w\,;\\
C(y,\,0)&=&C(0,\,w)=0.\\
\end{eqnarray*}
\item For every $y_1,\,y_2,\,w_1,\,w_2\in \mathbf{I}$ such that
$y_1\leq y_2$ and $w_1\leq w_2$,
$$C(y_2,\,w_2)-C(y_2,\,w_1)-C(y_1,\,w_2)+C(y_1,\,w_1)\geq 0,$$
\end{enumerate}
\end{defn}
\noindent that is, the C-measure of the box vertices lying in $\mathbf{I}^2$ is nonnegative. In particular, the last inequality implies  that $C(y,w)$ is increasing in both variables.
\\

Copulae are useful for capturing the dependence structure of random distributions with arbitrary marginals. This statement is clarified by Sklar's theorem which is now cited for the bivariate case.

\begin{thm}
Let $F(x_1,\,x_2)$ be the joint cumulative distribution function of random variables $X_1$ and $X_2$ having continuous marginal distributions $F_1(x_1)$ and $F_2(x_2)$. Then, there exists a unique bivariate copula
$C:\mathbf{I}^2\mapsto \mathbf{I}$ such that
\begin{equation}
F(x_1,\,x_2)=C\left(F_1(x_1),\,F_2(x_2)\right)\label{Sklar}
\end{equation}
where $C(\cdot,\,\cdot)$ is a joint distribution function with uniform marginals.
Conversely, for any continuous distribution function $F_1(x_1)$ and $F_2(x_2)$ and any copula $C$, the function $F$ defined in equation (\ref{Sklar}) is a joint distribution function with marginal distributions $F_1$ and $F_2$.
\end{thm}
Sklar's theorem provides a  scheme for constructing copulae. Indeed, the function
\begin{equation}
C(u_1,\,u_2)=F\left(F_1^{-1}(u_1),\,F_2^{-1}(u_2)\right)\label{C'=F}
\end{equation}
is a bivariate copula, where the quasi-inverse $F^{-1}_i$ for $i=1,\,2$ is defined by
\begin{equation}
F_i^{-1}(u)=\inf\{x|F_i(x)\geq u\}\quad\forall\,u\in (0,\,1)\label{F-1}.
\end{equation}
Much of the usefulness of the copulae follows from the fact that they are invariant with respect to strictly increasing  transformations.
More formally, let $X_1$ and $X_2$ be two continuous random variables with associated copula $C$. Now, letting $\alpha$ and $\beta$ be two strictly increasing functions and denoting by $C_{\alpha,\beta}$ the copula generated by $\alpha(X_1)$ and $\beta(X_2)$, it can be shown that for all $(u_1,\,u_2)\in \mathbf{I}^2,$
\begin{equation}
C_{\alpha,\,\beta}(u_1,u_2)=C(u_1,u_2).\label{C_transformation}
\end{equation}

Finally, let us denote by $c(\cdot,\cdot)$ the density function corresponding to the copula $C(\cdot,\cdot)$, that is,
\begin{equation*}
c(u_1,\,u_2)=\frac{\partial^2}{\partial u_1\partial u_2}C(u_1,u_2).
\end{equation*}
The following relationship between the joint density  $f(\cdot,\cdot)$  and  the copula density $c(\cdot,\cdot)$ can easily be obtained from equation (\ref{Sklar}):
\begin{equation}
f(x_1,x_2)=f_1(x_1)f_2(x_2)c(F_1(x_1),F_2(x_2))
\label{f=c}
\end{equation}
where $f_1(x_1)$ and $f_2(x_2)$ respectively denote the marginal density functions of $X_1$ and $X_2$. Thus, the copula density function can be 
expressed as follows:
\begin{equation}
c(u_1,\,u_2)=\frac{f(F_1^{-1}(u_1),\,F_2^{-1}(u_2))}{f_1(F_1^{-1}(u_1))\, f_2(F_2^{-1}(u_2))}.
\end{equation}

 Jaworski and Krzywda (2013) and Bosc (2012)  determined the copulae corresponding to certain correlated Brownian motions. Lager{\aa}s (2010) provides an explicit representation  of the copula associated with Brownian motion processes that are reflected at 0 and 1.  Several recent articles point out the usefulness of correlated Brownian motions and promote the use of copulae generated therefrom in connection with various applications. For instance, Chen et al. (2016) point out that correlated Brownian motions and their associated copulae can be utilized in the case of correlated assets occurring in risk management, pairs trading and multi-assets derivative's pricing. Deschatre (2016a,b) proposes to make use of asymmetric copulae generated from a Brownian motion and its reflection to model and control the distribution of their difference with applications to the energy market and the pricing of spread options.  

This paper which is principally based on the thesis of Vachon (2008),  is organized as follows. Several joint distributions related to certain  $\cal{BM}$  processes and their maxima are derived in the second section. The copulae associated with these joint distributions are then constructed  in the third section. 

\section{Brownian motion and related joint distributions}
As previously defined,  $\{W_t\}_{t\geq 0}$ shall denote a standard  $\cal{BM}$  and $M_t=\max_{0\leq s\leq t} W_s$, its maximum on the interval $[0,\,t]$. It is well known (see for instance, Etheridge (2002), Harrison (1990), Karlin and Taylor (1975), Revuz and Yor (2005), Rogers and Williams (2000)) that the joint distribution of $\left(W_t,\,M_t\right)$ and the marginal distribution of $M_t$ are respectively given by
\begin{equation}
\mathbb{P}\{M_t\leq a, W_t\leq x\}=\left\{\begin{array}{ll}
\Phi\left(\frac{x}{\sqrt{t}}\right)-\Phi\left(\frac{x-2a}{\sqrt{t}}\right) & \textrm{if }x\leq a\label{FC_MBS}\\
2\Phi\left(\frac{a}{\sqrt{t}}\right)-1 &\textrm{if } x>a
\end{array}\right.\quad\forall\,t\in\mathbb{R}_+
\end{equation}
and
\begin{equation}
\mathbb{P}\{M_t\leq a\}=2\Phi\left(\frac{a}{\sqrt{t}}\right)-1\quad\forall t\in \mathbb{R}_+,\label{F_Max}
\end{equation}
where $\Phi(\cdot)$ is the standard normal distribution function.

The first proposition of this section provides the joint distribution of $\{W_t^{(\mu,\,\sigma)}\}_{t\geq 0}$, a  $\cal{BM}$  with drift $\mu$ and variance $\sigma^2$, and $M_t^{(\mu,\,\sigma)}$, its maximum over the interval ${0\leq s\leq t}$.
This section conveniently provides detailed proofs of the distributional results  stated in the propositions, whereupon the corresponding copulae will be derived in the next section.

\begin{prop}\label{FC_MAX}\textit{(Harrison (1990))}
\begin{equation}
\mathbb{P}\{W^{(\mu,\sigma)}_t\leq x, M^{(\mu,\sigma)}_t\leq y\}=\left\{\begin{array}{ll}
\Phi\left(\frac{x-\mu t}{\sigma\sqrt{t}}\right)-e^{\frac{2\mu y}{\sigma^2}}\Phi\left(\frac{x-2y-\mu t}{\sigma\sqrt{t}}\right) & \textrm{\rm if $x\leq y$}\nonumber\\
\Phi\left(\frac{y-\mu t}{\sigma\sqrt{t}}\right)-e^{\frac{2\mu y}{\sigma^2}}\Phi\left(\frac{-y-\mu t}{\sigma\sqrt{t}}\right) & \textrm{\rm if $x>y$}.
\end{array}\right.
\end{equation}
\end{prop}

\ 

\noindent\begin{dem}  
	We first consider  the case where $x\leq y$. In light of equation (\ref{FC_MBS}), we have
	\begin{eqnarray*}
		\mathbb{P}\{W_t\in dx,M_t\leq y\} &=&\frac{1}{\sqrt{t}}\left(\phi\left(\frac{x}{\sqrt{t}}\right)-\phi\left(\frac{x-2y}{\sqrt{t}}\right)\right){\rm{d}}x,
	\end{eqnarray*}
	where $\phi(\cdot)$ denotes the standard normal density function.
	\\
	
	Now define
	$$\mathbb{Q}(A)=\int_A L_t(\omega){\rm d}\mathbb{P}(\omega),\quad A\in \mathfrak{F}_t^0$$
	where $L_t=e^{\mu W_t-\frac{1}{2}\mu^2 t}$ is the Radon-Nikodym's derivative of $\mathbb{Q}$ with respect to $\mathbb{P}$ and $\mathfrak{F}_t^0=\sigma(\{W_s, 0\leq s \leq t\})$ for all $t\in \mathbb{R}_+$, is the smallest $\sigma$-algebra generated by the ${\cal{BM}}$ up to time $t$. It follows from Girsanov's theorem that $\{W_t\}_{t\geq 0}$ is a  $\cal{BM}$  with drift $\mu$ under the new measure $\mathbb{Q}$.
	Therefore
	\begin{eqnarray*}
		\mathbb{Q}\{W_t\leq x, M_t\leq y\}&=&\int_{\{W_t\leq x,M_t\leq y\}}L_t(\omega){\rm d}\mathbb{P}(\omega)\nonumber\\
		&=&E^{\mathbb{P}}[\mathbf{1}_{\{W_t\leq x,M_t\leq y\}} L_t]\nonumber\\
		&=&\int_{-\infty}^x e^{\mu z-\frac{\mu^2 t}{2}} \frac{1}{\sqrt{t}}\left(\phi\left(\frac{z}{\sqrt{t}}\right)-\phi\left(\frac{z-2y}{\sqrt{t}}\right)\right){\rm{d}}z\nonumber\\
		&=&\int_{-\infty}^x \frac{1}{\sqrt{2\pi t}}e^{-\frac{(z-\mu t)^2}{2t}}{\rm{d}}z-e^{2\mu y}\int_{-\infty}^x \frac{1}{\sqrt{2\pi t}}e^{-\frac{(z-(2y+\mu t))^2}{2t}}{\rm{d}}z\nonumber\\
		&=&\Phi\left(\frac{x-\mu t}{\sqrt{t}}\right)-e^{2\mu y}\Phi\left(\frac{x-2y-\mu t}{\sqrt{t}}\right).\label{MAX2}
	\end{eqnarray*}
\end{dem}

\
\begin{rem} Note that the marginal distribution of $M_t$, which is given by
\begin{equation}
\mathbb{P}\{M_t^{(\mu,\sigma)}\leq y\}=\Phi\left(\frac{y-\mu t}{\sigma\sqrt{t}}\right)-e^{\frac{2\mu y}{\sigma^2}}\Phi\left(\frac{-y-\mu t}{\sigma\sqrt{t}}\right),\label{F_Max2}
\end{equation}
can easily be derived from Proposition~\ref{FC_MAX} since for $x>y$, $\{M_t^{(\mu,\sigma)}\leq y\}\subset\{W_t^{(\mu,\sigma)}\leq x\}$.
\end{rem}
One can generalize these last results by making use of the following properties of the multivariate normal distribution:
\begin{equation}
\Phi_2(z_1,\,z_2;\,\rho)=\Phi_2(z_2,\,z_1;\,\rho),\label{Phi2_1}
\end{equation}
\begin{equation}
\Phi(z_1)-\Phi_2(z_1,\,z_2;\,\rho)=\Phi_2(z_1,\,-z_2;\,-\rho),\label{Phi2_2}
\end{equation}
\begin{eqnarray}
\Phi_3(z_1,\,z_2,\,z_3;\,\rho_{12},\,\rho_{13},\,\rho_{23})&=&\Phi_3(z_2,\,z_1,\,z_3;\,\rho_{12},\,\rho_{23},\,\rho_{13})\label{Phi3_1}\\
&=&\Phi_3(z_3,\,z_1,\,z_2;\,\rho_{13},\,\rho_{23},\,\rho_{12})\nonumber
\end{eqnarray}
and
\begin{eqnarray}
\lefteqn{\Phi_2(z_2,\,z_3;\,\rho_{23})-\Phi_3(z_1,\,z_2,\,z_3;\,\rho_{12},\,\rho_{13},\,\rho_{23})}\label{Phi_2_3}\\
&=&\Phi_3(-z_1,\,z_2,\,z_3;\,-\rho_{12},\,-\rho_{13},\,\rho_{23}).\nonumber
\end{eqnarray}

\begin{Lem}\label{Phi_2_prop} Let $z_1,\,z_2,\,z_3$ be real constants and $\rho\geq 0$. If $z_1=-\rho z_2+\sqrt{1-\rho^2}\,z_3$, then
\begin{equation}
\Phi_2(z_1,\,z_2;\,-\rho)+\Phi_2(-z_1,\,z_3;\,-\sqrt{1-\rho^2})=\Phi(z_2)\Phi(z_3)\label{Phi2_3}
\end{equation}
and
\begin{equation}
\Phi_2(z_1,\,z_2;\,-\rho)+\Phi(-z_2)\Phi(z_3)=\Phi_2(z_1,\,z_3;\,\sqrt{1-\rho^2}).\label{Phi2_4}
\end{equation}
\end{Lem}

\noindent \begin{dem}  Let $Z_2$ and $Z_3$ be two independent standard normal random variables, and $Z_1$ a random variable defined by $Z_1=-\rho Z_2+\sqrt{1-\rho^2}Z_3$. Note that $Z_1$ has also a standard normal distribution, the random vectors $(Z_1,\,Z_2)$ and $(-Z_1,\,Z_3)$ have bivariate normal distributions with correlation coefficients given by $-\rho$ and $-\sqrt{1-\rho^2}$, respectively. Then,
\begin{eqnarray*}
\lefteqn{\Phi_2(z_1,\,z_2;\,-\rho)+\Phi_2(-z_1,\,z_3;\,-\sqrt{1-\rho^2})}\\
&=&\mathbb{P}\{Z_1\leq z_1,\,Z_2\leq z_2\}+\mathbb{P}\{-Z_1\leq -z_1,\,Z_3\leq z_3\}\\
&=&\mathbb{P}\{Z_1\leq z_1,\,Z_2\leq z_2,\,Z_3\leq z_3\}+\mathbb{P}\{Z_1\leq z_1,\,Z_2\leq z_2,\,Z_3\geq z_3\}{}\\
& &{}+\mathbb{P}\{Z_1\geq z_1,\,Z_2\leq z_2,\,Z_3\leq z_3\}+\mathbb{P}\{Z_1\geq z_1,\,Z_2\geq z_2,\,Z_3\leq z_3\}.
\end{eqnarray*}
We now replace $Z_1$ by $-\rho Z_2+\sqrt{1-\rho^2}Z_3$ and $z_1$ by $-\rho z_2+\sqrt{1-\rho}z_3$. Since the events $\{Z_1\leq z_1,\,Z_2\leq z_2,\,Z_3\geq z_3\}$ and $\{Z_1\geq z_1,\,Z_2\geq z_2,\,Z_3\leq z_3\}$ are clearly empty, we obtain
\begin{eqnarray*}
\lefteqn{\Phi_2(z_1,\,z_2;\,-\rho)+\Phi_2(-z_1,\,z_3;\,-\sqrt{1-\rho^2})}\\
&=&\mathbb{P}\{Z_1\leq z_1,\,Z_2\leq z_2,\,Z_3\leq z_3\}+\mathbb{P}\{Z_1\geq z_1,\,Z_2\leq z_2,\,Z_3\leq z_3\}\\
&=&\mathbb{P}\{Z_2\leq z_2,\,Z_3\leq z_3\}=\Phi(z_2)\Phi(z_3).
\end{eqnarray*}
It follows from equations $(\ref{Phi2_2})$ and $(\ref{Phi2_3})$ that
\begin{eqnarray*}
\lefteqn{\Phi_2(z_1,\,z_2;\,-\rho)+\Phi_2(-z_1,\,z_3;\,-\sqrt{1-\rho^2})=\Phi(z_2)\Phi(z_3)}\\
&\Rightarrow&\Phi_2(z_1,\,z_2;\,-\rho)+\Phi(z_3)-\Phi_2(z_1,\,z_3;\,\sqrt{1-\rho^2})=(1-\Phi(-z_2))\Phi(z_3)\\
&\Rightarrow&\Phi_2(z_1,\,z_2;\,-\rho)+\Phi(-z_2)\Phi(z_3)=\Phi_2(z_1,\,z_3;\,\sqrt{1-\rho^2}).
\end{eqnarray*}
\end{dem}\\
\\
The joint distributions that will be  considered  further involve  integrals for which closed form representations are given in the next proposition. 
\begin{prop}\label{Int_Phi3}Let $a$, $h$, $\theta_i$, $i=1,\,2,\,3$, $\delta_j$ and $\eta_j >0$, $j=0,\,1,\,2,\,3$ be constant, and $\mathbf{R}=[\rho_{ij}]_{i,\,j=1,\,2,\,3}$ be a correlation  matrix, then
\begin{eqnarray}
\lefteqn{\int_{-\infty}^a\exp{(hs)}\Phi_3\left(\frac{\delta_1+\theta_1s}{\eta_1},\,\frac{\delta_2+\theta_2s}{\eta_2},\,\frac{\delta_3+\theta_3s}{\eta_3};\mathbf{R}\right)\phi\left(\frac{s-\delta_0}{\eta_0}\right)\frac{{\rm{d}}s}{\eta_0}}\label{Phi4_1}\\
&=& \exp{\left(h\delta_0+\frac{h^2\eta_0^2}{2}\right)}\Phi_4\left(\frac{a-\delta_0^*}{\eta_0},\,\frac{\delta_1+\theta_1\delta_0^*}{\kappa_1},\,\frac{\delta_2+\theta_2\delta_0^*}{\kappa_2},\,\frac{\delta_3+\theta_3\delta_0^*}{\kappa_3};\mathbf{R}^*\right)\nonumber
\end{eqnarray}
and
\begin{eqnarray}
\lefteqn{\int_a^{+\infty}\exp{(hs)}\Phi_3\left(\frac{\delta_1+\theta_1s}{\eta_1},\,\frac{\delta_2+\theta_2s}{\eta_2},\,\frac{\delta_3+\theta_3s}{\eta_3};\mathbf{R}\right)\phi\left(\frac{s-\delta_0}{\eta_0}\right)\frac{{\rm{d}}s}{\eta_0}}\label{Phi4_2}\\
&=& \exp{\left(h\delta_0+\frac{h^2\eta_0^2}{2}\right)}\Phi_4\left(\frac{-a+\delta_0^*}{\eta_0},\,\frac{\delta_1+\theta_1\delta_0^*}{\kappa_1},\,\frac{\delta_2+\theta_2\delta_0^*}{\kappa_2},\,\frac{\delta_3+\theta_3\delta_0^*}{\kappa_3};\mathbf{R}^{**}\right)\nonumber
\end{eqnarray}
where $\delta_0^*=\delta_0+h\eta_0^2$; $\kappa_i=\sqrt{\theta_i^2\eta_0^2+\eta_i^2}$ for $i=1,\,2,\,3$;  $\mathbf{R}^*=[\rho_{ij}^*]_{i,\,j=1,\,2,\,3,\,4}$ with $\rho_{1\,i+1}^*=-{\theta_i\eta_0}/{\kappa_i}$,  $i=1,\,2,\,3$; $\rho_{2\,i+1}^*=({\rho_{1i}\eta_1\eta_i+\theta_1\theta_i\eta_0^2})/({\kappa_1\kappa_i})$, $i=2,\,3$; $\rho^*_{34}=({\rho_{23}\eta_2\eta_3+\theta_2\theta_3\eta_0^2})/({\kappa_2\kappa_3})$; and finally $\mathbf{R}^{**}=[\rho_{ij}^{**}]_{i,\,j=1,\,2,\,3,\,4}$ with $\rho_{1i}^{**}=-\rho_{1i}^*$, $i=2,\,3,\,4$;  $\rho_{ij}^{**}=\rho_{ij}^*$, $i,j=2,3,4$.
\end{prop}

These results are established by making use of properties of the conditional multivariate normal distribution. Note that this proposition is related to a result appearing in Lee (2003) whose derivation relies on the {Esscher transform}.

\

\noindent \begin{dem}  Let $\mathbf{X}=(X_1,\,X_2,\,X_3,\,X_4)'$ be a normally distributed random vector  such $E[X_i]=\mu_i$, $Var[X_i]=\sigma_i^2$ and $\mathbf{R}^*=[\rho_{ij}^*]$ for $i,\,j=1,\,2,\,3,\,4$. Then the conditional distribution of $(X_2,\,X_3,\,X_4)$ given $X_1=x_1$ is a trivariate normal distribution (Anderson, 2003) with mean vector 
\begin{eqnarray*}
\mathbf{\mu}^{(1)}+\mathbf{\Sigma}_{12}\mathbf{\Sigma}_{22}^{-1}\left(\mathbf{x}^{(2)}
-\mathbf{\mu}^{(2)}\right)&=&\left(\begin{array}{c}
\mu_2+\frac{\sigma_2}{\sigma_1}\rho_{12}\left(x_1-\mu_1\right)\\
\mu_3+\frac{\sigma_3}{\sigma_1}\rho_{13}\left(x_1-\mu_1\right)\\
\mu_4+\frac{\sigma_4}{\sigma_1}\rho_{14}\left(x_1-\mu_1\right)
\end{array}\right)
\end{eqnarray*}
and covariance matrix 
\begin{eqnarray*}
\mathbf{\Sigma}_{11}-\mathbf{\Sigma}_{12}\mathbf{\Sigma}_{22}^{-1}\mathbf{\Sigma}_{21}
&=&\left(\begin{array}{ccc}
\sigma_2^2(1-\rho_{12}^2) & \sigma_2\sigma_3(\rho_{23}-\rho_{12}\rho_{13}) & \sigma_2\sigma_4(\rho_{24}-\rho_{12}\rho_{14})\\
\sigma_2\sigma_3(\rho_{23}-\rho_{12}\rho_{13}) & \sigma_3^2(1-\rho_{13}^2) & \sigma_3\sigma_4(\rho_{34}-\rho_{13}\rho_{14})\\
\sigma_2\sigma_4(\rho_{24}-\rho_{12}\rho_{14}) &\sigma_3\sigma_4(\rho_{34}-\rho_{13}\rho_{14})  & \sigma_4^2(1-\rho_{14}^2)
\end{array}\right).
\end{eqnarray*}

Thus,
\begin{eqnarray}
\lefteqn{\Phi_4\left(\frac{x_1-\mu_1}{\sigma_1},\,\frac{x_2-\mu_2}{\sigma_2},\,\frac{x_3-\mu_3}{\sigma_3},\,\frac{x_4-\mu_4}{\sigma_4};\mathbf{R}^*\right)}\nonumber\\
&=&\mathbb{P}\{X_1\leq x_1,\,X_2\leq x_2,\,X_3\leq x_3,\,X_4\leq x_4\}\nonumber\\
&=&\int_{-\infty}^{x_1} \mathbb{P}\{X_2\leq x_2,\,X_3\leq x_3,\,X_4\leq x_4\,|\,X_1=s\}\mathbb{P}\left\{X_1\in {\rm{d}}s\right\}\nonumber\\
&=&\int_{-\infty}^{x_1}\Phi_3\left(\frac{x_2-(\mu_2+\rho_{12}^*\frac{\sigma_2}{\sigma_1}(s-\mu_1))}{\sigma_2\sqrt{1-(\rho_{12}^*)^{2}}},\,\ldots,\,\frac{x_4-(\mu_4+\rho_{14}^*\frac{\sigma_4}{\sigma_1}(s-\mu_1))}{\sigma_4\sqrt{1-(\rho_{14}^*)^2}};\right.{}\label{Phi4_3}\\ & &{}\left.\frac{\rho_{23}^*-\rho_{12}^*\rho_{13}^*}{\sqrt{1-(\rho_{12}^*)^2}\sqrt{1-(\rho_{13}^*)^2}},\,\frac{\rho_{24}^*-\rho_{12}^*\rho_{14}^*}{\sqrt{1-(\rho_{12}^*)^2}\sqrt{1-(\rho_{14}^*)^2}},\,\frac{\rho_{34}^*-\rho_{13}^*\rho_{14}^*}{\sqrt{1-(\rho_{13}^*)^2}\sqrt{1-(\rho_{14}^*)^2}}\right){}\nonumber\\
& &{}\times\phi\left(\frac{s-\mu_1}{\sigma_1}\right)\frac{{\rm{d}}s}{\sigma_1}.\nonumber
\end{eqnarray}
Now letting $x_1=a$, $x_{i+1}=\delta_{i}$, $\mu_1=\delta_0$, $\mu_{i+1}=-\theta_{i}\delta_0$, $\sigma_1=\eta_0$, $\sigma_{i+1}=\kappa_{i}$ for $i=1,\,2,\,3$ and replacing in equation (\ref{Phi4_3}) the elements of the matrix $\mathbf{R}^*$ with their respective values, we obtain
\begin{eqnarray*} \lefteqn{\Phi_4\left(\frac{a-\delta_0}{\eta_0},\,\frac{\delta_1+\theta_1\delta_0}{\kappa_1},\,\frac{\delta_2+\theta_2\delta_0}{\kappa_2},\,\frac{\delta_3+\theta_3\delta_0}{\kappa_3};\mathbf{R}^{*}\right)}\\
&=&\int_{-\infty}^a\Phi_3\left(\frac{\delta_1+\theta_1s}{\eta_1},\,\frac{\delta_2+\theta_2s}{\eta_2},\,\frac{\delta_3+\theta_3s}{\eta_3};\mathbf{R}\right)\phi\left(\frac{s-\delta_0}{\eta_0}\right)\frac{{\rm{d}}s}{\eta_0}.
\end{eqnarray*}
This establishes equation  (\ref{Phi4_1}) for $h=0$.

The case $h\neq 0$ follows from the last expression by completing the square  in the exponent of  $\exp{(hs)}\phi\left(\frac{s-\delta_0}{\eta_0}\right)$, so that
$$\exp{(hs)}\phi\left(\frac{s-\delta_0}{\eta_0}\right)=\exp{(h\delta_0+\frac{h^2\eta_0^2}{2})}\,\phi\left(\frac{s-\delta_0^*}{\eta_0}\right).$$

Finally, the last result, that is, equation (\ref{Phi4_2}) is similarly obtained  on noting that
\begin{eqnarray*}
\lefteqn{\Phi_4\left(\frac{-x_1+\mu_1}{\sigma_1},\,\frac{x_2-\mu_2}{\sigma_2},\,\frac{x_3-\mu_3}{\sigma_3},\,\frac{x_4-\mu_4}{\sigma_4};\mathbf{R}^{**}\right)}\\
&=&\mathbb{P}\{(-X_1)\leq -x_1,\,X_2\leq x_2,\,X_3\leq x_3,\,X_4\leq x_4\}\\
&=&\mathbb{P}\{X_1\geq x_1,\,X_2\leq x_2,\, X_3\leq x_3,\,X_4\leq x_4\}.
\end{eqnarray*}
\end{dem}

Additionally, as $\delta_3\rightarrow\infty$, we have
\begin{eqnarray}
\lefteqn{\int_{-\infty}^a\exp{(hs)}\Phi_2\left(\frac{\delta_1+\theta_1s}{\eta_1},\,\frac{\delta_2+\theta_2s}{\eta_2};\,\rho_{12}\right)}\label{Int2_1}\\
& &{}\times\phi\left(\frac{s-\delta_0}{\eta_0}\right)\frac{{\rm{d}}s}{\eta_0}\nonumber\\
&=& \exp{\left(h\delta_0+\frac{h^2\eta_0^2}{2}\right)}\Phi_3\left(\frac{a-\delta_0^*}{\eta_0},\,\frac{\delta_1+\theta_1\delta_0^*}{\kappa_1};\,\rho_{12}^{*},\,\rho_{13}^*,\,\rho_{23}^*\right)\nonumber
\end{eqnarray}
and
\begin{eqnarray}
\lefteqn{\int_a^{+\infty}\exp{(hs)}\Phi_2\left(\frac{\delta_1+\theta_1s}{\eta_1},\,\frac{\delta_2+\theta_2s}{\eta_2};\,\rho_{12}\right)}\label{Int2_2}\\
& &{}\times\phi\left(\frac{s-\delta_0}{\eta_0}\right)\frac{{\rm{d}}s}{\eta_0}\nonumber\\
&=& \exp{\left(h\delta_0+\frac{h^2\eta_0^2}{2}\right)}\Phi_3\left(\frac{-a+\delta_0^*}{\eta_0},\,\frac{\delta_1+\theta_1\delta_0^*}{\kappa_1};\,-\rho_{12}^{*},\,-\rho_{13}^*,\,\rho_{23}^*\right).\nonumber
\end{eqnarray}
Similarly, as $\delta_2\rightarrow\infty$, it follows from equations  (\ref{Int2_1}) and (\ref{Int2_2}) that
\begin{eqnarray}
\lefteqn{\int_{-\infty}^a\exp{(hs)}\Phi\left(\frac{\delta_1+\theta_1s}{\eta_1}\right)\phi\left(\frac{s-\delta}{\eta}\right)\frac{{\rm{d}}s}{\eta}}\label{Int1_1}\\
&=& \exp{\left(h\delta+\frac{h^2\eta^2}{2}\right)}\Phi_2\left(\frac{a-\delta^*}{\eta},\,\frac{\delta_1+\theta_1\delta^*}{\kappa_1};\,-\frac{\theta_1\eta}{\kappa_1}\right)\nonumber
\end{eqnarray}
and
\begin{eqnarray}
\lefteqn{\int_a^{+\infty}\exp{(hs)}\Phi\left(\frac{\delta_1+\theta_1s}{\eta_1}\right)\phi\left(\frac{s-\delta}{\eta}\right)\frac{{\rm{d}}s}{\eta}}\label{Int1_2}\\
&=& \exp{\left(h\delta+\frac{h^2\eta^2}{2}\right)}\Phi_2\left(\frac{-a+\delta^*}{\eta},\,\frac{\delta_1+\theta_1\delta^*}{\kappa_1};\,\frac{\theta_1\eta}{\kappa_1}\right).\nonumber
\end{eqnarray}

These results enable one to establish 
the distribution of $\left(W_T^{(\mu,\,\sigma)}, \,M_t^{(\mu,\,\sigma)}\right)$ within the interval  $0<t\leq T$ as specified in the next proposition.
\begin{prop}\label{M_t_W_T} \textit{(Chuang (1996) and Lee (2003))}
\begin{eqnarray}
\mathbb{P}\{W_T^{(\mu,\sigma)}\leq x,\,M_t^{(\mu,\,\sigma)}\leq y\}&=&\Phi_2\left(\frac{x-\mu T}{\sigma\sqrt{T}},\,\frac{y-\mu t}{\sigma\sqrt{t}};\,\sqrt{\frac{t}{T}}\right){}\\
& &{}-e^{\frac{2\mu y}{\sigma^2}}\Phi_2\left(\frac{x-2y-\mu T}{\sigma\sqrt{T}},\,\frac{-y-\mu t}{\sigma\sqrt{t}};\,\sqrt{\frac{t}{T}}\right).\nonumber
\end{eqnarray}
\end{prop}

\

\noindent\begin{dem}
	\begin{eqnarray}
	\lefteqn{P\{W_T^{(\mu,\sigma)}\leq x,\,M_t^{(\mu,\sigma)}\leq y\}}\nonumber\\
	&=&\int_{-\infty}^{+\infty}P\{W_T^{(\mu,\sigma)}\leq x,\,M_t^{(\mu,\sigma)}\leq y|W_T^{(\mu,\,\,\sigma)}-W^{(\mu,\sigma)}_t=z\}{}\nonumber\\
	& &{}\times P\{W_T^{(\mu,\,\,\sigma)}-W^{(\mu,\sigma)}_t\in {\rm{d}}z\}\nonumber\\
	&=&\int_{-\infty}^{+\infty}P\{W_t^{(\mu,\sigma)}\leq x-z,\,M_t^{(\mu,\sigma)}\leq y|W_T^{(\mu,\,\,\sigma)}-W^{(\mu,\sigma)}_t=z\}{}\nonumber\\
	& &{}\times \phi\left(\frac{z-\mu(T-t)}{\sigma\sqrt{T-t}}\right)\frac{{\rm{d}}z}{\sigma\sqrt{T-t}}\nonumber\\
	&=&\int_{-\infty}^{+\infty}P\{W_t^{(\mu,\sigma)}\leq x-z,\,M_t^{(\mu,\sigma)}\leq y\}\phi\left(\frac{z-\mu(T-t)}{\sigma\sqrt{T-t}}\right)\frac{{\rm{d}}z}{\sigma\sqrt{T-t}}\nonumber\\
	&=&\int_{-\infty}^{x-y}P\{W_t^{(\mu,\sigma)}\leq x-z,\,M_t^{(\mu,\sigma)}\leq y\}\phi\left(\frac{z-\mu(T-t)}{\sigma\sqrt{T-t}}\right)\frac{{\rm{d}}z}{\sigma\sqrt{T-t}}{}\label{Mt_1}\\
	& &{}+\int_{x-y}^{+\infty}P\{W_t^{(\mu,\sigma)}\leq x-z,\,M_t^{(\mu,\sigma)}\leq y\}\phi\left(\frac{z-\mu(T-t)}{\sigma\sqrt{T-t}}\right)\frac{{\rm{d}}z}{\sigma\sqrt{T-t}}.\nonumber
	\end{eqnarray}
	By replacing in equation (\ref{Mt_1}) the result of Proposition~\ref{FC_MAX}, we obtain for the first part of the equation:
	\begin{eqnarray}
	\lefteqn{\int_{-\infty}^{x-y}P\{W_t^{(\mu,\sigma)}\leq x-z,\,M_t^{(\mu,\sigma)}\leq y\}\phi\left(\frac{z-\mu(T-t)}{\sigma\sqrt{T-t}}\right)\frac{{\rm{d}}z}{\sigma\sqrt{T-t}}}\nonumber\\
	&=&\Phi\left(\frac{y-\mu t}{\sigma\sqrt{t}}\right)\Phi\left(\frac{x-y-\mu(T-t)}{\sigma\sqrt{T-t}}\right){}\label{Mt_2}\\
	& &{}-e^{\frac{2\mu y}{\sigma^2}}\Phi\left(\frac{-y-\mu t}{\sigma\sqrt{t}}\right)\Phi\left(\frac{x-y-\mu(T-t)}{\sigma\sqrt{T-t}}\right);\nonumber
	\end{eqnarray}
	as for the second part,
	\begin{eqnarray}
	\lefteqn{\int_{x-y}^{+\infty}P\{W_t^{(\mu,\sigma)}\leq x-z,\,M_t^{(\mu,\sigma)}\leq y\}\phi\left(\frac{z-\mu(T-t)}{\sigma\sqrt{T-t}}\right)\frac{{\rm{d}}z}{\sigma\sqrt{T-t}}}\nonumber\\
	&=&\int_{x-y}^{+\infty}\Phi\left(\frac{-z-(\mu t-x)}{\sigma\sqrt{t}}\right)\phi\left(\frac{z-\mu(T-t)}{\sigma\sqrt{T-t}}\right)\frac{{\rm{d}}z}{\sigma\sqrt{T-t}}{}\nonumber\\
	& &{}-e^{\frac{2\mu y}{\sigma^2}}\int_{x-y}^{+\infty}\Phi\left(\frac{-z-(2y+\mu t-x)}{\sigma\sqrt{t}}\right)\phi\left(\frac{z-\mu(T-t)}{\sigma\sqrt{T-t}}\right)\frac{{\rm{d}}z}{\sigma\sqrt{T-t}}\nonumber\\
	&=&\Phi_2\left(\frac{y-x+\mu(T-t)}{\sigma\sqrt{T-t}},\,\frac{x-\mu T}{\sigma\sqrt{T}};\,-\sqrt{1-\frac{t}{T}}\right){}\label{Mt_3}\\
	& &-e^{\frac{2\mu y}{\sigma^2}}\Phi_2\left(\frac{y-x+\mu(T-t)}{\sigma\sqrt{T-t}},\,\frac{x-2y-\mu T}{\sigma\sqrt{T}};\,-\sqrt{1-\frac{t}{T}}\right)\nonumber
	\end{eqnarray}
	where the last equality follows from equation (\ref{Int1_2}). On combining the  last two results and applying Lemma~\ref{Phi_2_prop}, we obtain
	\begin{eqnarray*}
		\lefteqn{P\{W_T^{(\mu,\sigma)}\leq x,\,M_t^{(\mu,\sigma)}\leq y\}}\nonumber\\
		&=&\Phi_2\left(\frac{x-\mu T}{\sigma\sqrt{T}},\,\frac{y-\mu t}{\sigma\sqrt{t}};\,\sqrt{\frac{t}{T}}\right)-e^{\frac{2\mu y}{\sigma^2}}\Phi_2\left(\frac{x-2y-\mu T}{\sigma\sqrt{T}},\,\frac{-y-\mu t}{\sigma\sqrt{t}};\,\sqrt{\frac{t}{T}}\right).
	\end{eqnarray*}
\end{dem}\\
\

\begin{rem} As expected, when $t\rightarrow T$, $\mathbb{P}\{W_T^{(\mu,\,\sigma)}\leq x,\,M_t^{(\mu,\,\sigma)}\leq y\}\rightarrow\mathbb{P}\{W_T^{(\mu,\,\sigma)}\leq x,\,M_T^{(\mu,\,\sigma)}\leq y\}$.
\end{rem}
\
Next, the joint distribution of $W_T^{(\mu,\,\sigma)}$ and $M_{(s,\,t)}^{(\mu,\,\sigma)}$, where $M_{(s,t)}^{(\mu,\,\sigma)}=\max_{s\leq u\leq t}W_u^{(\mu,\,\sigma)}$ and $0<s<t\leq T$, is considered.

\begin{prop}\label{M_st}\textit{(Lee (2003))}
\begin{eqnarray}
\lefteqn{\mathbb{P}\{W_T^{(\mu,\sigma)}\leq x,\,M_{(s,\,t)}^{(\mu,\,\sigma)}\leq y\}}\nonumber\\
&=&\Phi_3\left(\frac{x-\mu T}{\sigma\sqrt{T}},\,\frac{y-\mu t}{\sigma\sqrt{t}},\,\frac{y-\mu s}{\sigma\sqrt{s}};\,\sqrt{\frac{t}{T}},\,\sqrt{\frac{s}{T}},\,\sqrt{\frac{s}{t}}\right){}\label{Mst_2}\\
& &{}-e^{\frac{2\mu y}{\sigma^2}}\Phi_3\left(\frac{x-2y-\mu T}{\sigma\sqrt{T}},\,\frac{-y-\mu t}{\sigma\sqrt{t}},\,\frac{y+\mu s}{\sigma\sqrt{s}};\,\sqrt{\frac{t}{T}},\,-\sqrt{\frac{s}{T}},\,-\sqrt{\frac{s}{t}}\right)\nonumber
\end{eqnarray}
and
\begin{eqnarray}
\mathbb{P}\{M_{(s,\,t)}^{(\mu,\,\sigma)}\leq y\}&=&\Phi_2\left(\frac{y-\mu t}{\sigma\sqrt{t}},\,\frac{y-\mu s}{\sigma\sqrt{s}};\,\sqrt{\frac{s}{t}}\right){}\label{Mst_1}\\
& &{}-e^{\frac{2\mu y}{\sigma^2}}\Phi_2\left(\frac{-y-\mu t}{\sigma\sqrt{t}},\,\frac{y+\mu s}{\sigma\sqrt{s}};\,-\sqrt{\frac{s}{t}}\right).\nonumber
\end{eqnarray}
\end{prop}
\


\noindent\begin{dem}
	Let us first consider the joint distribution of a ${\cal{BM}}$ and its maximum on the interval $[s,\,t]$. In that case,
	\begin{eqnarray}
	\lefteqn{\mathbb{P}\{W_T^{(\mu,\sigma)}\leq x,\,M_{(s,\,t)}^{(\mu,\,\sigma)}\leq y\}}\nonumber\\
	&=&\int_{-\infty}^y\mathbb{P}\{W_T^{(\mu,\sigma)}\leq x,\,M_{(s,\,t)}^{(\mu,\,\sigma)}\leq y|W_s^{(\mu,\,\sigma)}=z\}\mathbb{P}\{W_s^{(\mu,\,\sigma)}\in {\rm{d}}z\}\nonumber\\
	&=&\int_{-\infty}^y\mathbb{P}\{W_T^{(\mu,\sigma)}-W_s^{(\mu,\sigma)}\leq x-z,\,M_{(s,\,t)}^{(\mu,\,\sigma)}-W_s^{(\mu,\sigma)}\leq y-z|W_s^{(\mu,\,\sigma)}=z\}{}\nonumber\\
	& &{}\times\phi\left(\frac{z-\mu s}{\sigma\sqrt{s}}\right)\frac{{\rm{d}}z}{\sigma\sqrt{s}}\nonumber\\
	&=&\int_{-\infty}^y\mathbb{P}\{W_T^{(\mu,\sigma)}-W_s^{(\mu,\sigma)}\leq x-z,\,\max_{s\leq u\leq t}\{W^{(\mu,\,\sigma)}_u-W_s^{(\mu,\sigma)}\}\leq y-z\}{}\nonumber\\
	& &{}\times\phi\left(\frac{z-\mu s}{\sigma\sqrt{s}}\right)\frac{{\rm{d}}z}{\sigma\sqrt{s}}\nonumber
	\end{eqnarray}
	\begin{eqnarray}
	&=&\int_{-\infty}^y\mathbb{P}\{W_{T-s}^{(\mu,\sigma)}\leq x-z,\,\max_{s\leq u\leq t}W^{(\mu,\,\sigma)}_{u-s}\leq y-z\}\phi\left(\frac{z-\mu s}{\sigma\sqrt{s}}\right)\frac{{\rm{d}}z}{\sigma\sqrt{s}}\nonumber\\
	&=&\int_{-\infty}^y\mathbb{P}\{W_{T-s}^{(\mu,\sigma)}\leq x-z,\,\max_{0\leq v\leq t-s}W^{(\mu,\,\sigma)}_v\leq y-z\}\phi\left(\frac{z-\mu s}{\sigma\sqrt{s}}\right)\frac{{\rm{d}}z}{\sigma\sqrt{s}}\nonumber\\
	&=&\int_{-\infty}^y\mathbb{P}\{W_{T-s}^{(\mu,\sigma)}\leq x-z,\,M^{(\mu,\,\sigma)}_{t-s}\leq y-z\}\phi\left(\frac{z-\mu s}{\sigma\sqrt{s}}\right)\frac{{\rm{d}}z}{\sigma\sqrt{s}}.\label{Mst_3a}
	\end{eqnarray}
	On applying the result of Proposition~\ref{M_t_W_T} to the first term in the integrand of (\ref{Mst_3a}), we obtain
	\begin{eqnarray*}
		\lefteqn{\mathbb{P}\{W_T^{(\mu,\sigma)}\leq x,\,M_{(s,\,t)}^{(\mu,\,\sigma)}\leq y\}}\nonumber\\
		&=&\int_{-\infty}^y\Phi_2\left(\frac{x-z-\mu (T-s)}{\sigma\sqrt{T-s}},\,\frac{y-z-\mu (t-s)}{\sigma\sqrt{t-s}};\,\sqrt{\frac{t-s}{T-s}}\right)\label{Mst_3}{}\\
		& &{}\times\phi\left(\frac{z-\mu s}{\sigma\sqrt{s}}\right)\frac{{\rm{d}}z}{\sigma\sqrt{s}}{}\nonumber\\
		& &{}-\int_{-\infty}^y e^{\frac{2\mu (y-z)}{\sigma^2}}\Phi_2\left(\frac{x+z-2y-\mu (T-s)}{\sigma\sqrt{T-s}},\,\frac{-y+z-\mu (t-s)}{\sigma\sqrt{t-s}};\,\sqrt{\frac{t-s}{T-s}}\right){}\nonumber\\
		& &{}\ \ \ \ \times\phi\left(\frac{z-\mu s}{\sigma\sqrt{s}}\right)\frac{{\rm{d}}z}{\sigma\sqrt{s}}\nonumber\\
		&=&\Phi_3\left(\frac{x-\mu T}{\sigma\sqrt{T}},\,\frac{y-\mu t}{\sigma\sqrt{t}},\,\frac{y-\mu s}{\sigma\sqrt{s}};\,\sqrt{\frac{t}{T}},\,\sqrt{\frac{s}{T}},\,\sqrt{\frac{s}{t}}\right)\label{Mst_4}{}\\
		& &{}-e^{\frac{2\mu y}{\sigma^2}}\Phi_3\left(\frac{x-2y-\mu T}{\sigma\sqrt{T}},\,\frac{-y-\mu t}{\sigma\sqrt{t}},\,\frac{y+\mu s}{\sigma\sqrt{s}};\,\sqrt{\frac{t}{T}},\,-\sqrt{\frac{s}{T}},\,-\sqrt{\frac{s}{t}}\right)\nonumber
	\end{eqnarray*}
	where the last equation follows from equation (\ref{Int2_1}).
	
	Finally, we obtain result (\ref{Mst_1}) by  letting $x$ tend to  $+ \infty $ in equation (\ref{Mst_2}).
\end{dem}

\

\begin{rem} When $s\rightarrow 0$, $\mathbb{P}\{W_T^{(\mu,\,\sigma)}\leq x,\,M_{(s,\,t)}^{(\mu,\,\sigma)}\leq y\}\rightarrow\mathbb{P}\{W_T^{(\mu,\,\sigma)}\leq x,\,M_t^{(\mu,\,\sigma)}\leq y\}$ and $\mathbb{P}\{M_{(s,\,t)}^{(\mu,\,\sigma)}\leq y\}\rightarrow\mathbb{P}\{M_t^{(\mu,\,\sigma)}\leq y\}$.
\end{rem}

\

Consider $\mathbf{W}=(W^1,\,W^2)'$ a Brownian vector where $W^1$ and $W^2$ are two independent standard ${\cal{BM}}$ processes. 
On letting $\{B^1_t=\sigma_1 (\rho W^2_t+\sqrt{1-\rho^2}W^1_t)+\mu_1 t\}_{t\in\mathbb{R}_+}$ and $\{B^2_t=\sigma_2 W^2_t+\mu_2 t\}_{t\in\mathbb{R}_+}$, one can construct a correlated two-dimensional  ${\cal{BM}}$ process. Then, it can easily be verified that $\{B^i_t\}_{t\in\mathbb{R}_+}$ is a $(\mu_i,\,\sigma_i)$-${\cal{BM}}$ for $i=1,\,2$ and the correlation between $B^1_t$ and $B^2_t$ is equal to $\rho$.
We say that $(B^1,\,B^2)'$ is a $\left(\mathbf{\mu},\,\mathbf{\Sigma}\right)$-${\cal{BM}}$ with drift vector
$$\mathbf{\mu}=(
\mu_1, \mu_2)'$$
and covariance matrix
$$\mathbf{\Sigma}=\left(\begin{array}{cc}
\sigma_1^2 & \rho\sigma_1\sigma_2\\
\rho\sigma_1\sigma_2 & \sigma_2^2
\end{array}\right).$$

Finally, we consider  the joint distribution of $\left(B^1_t,M^2_{(s,\,t)}\right)$ for correlated ${\cal{BM}}$s where $M^2_{(s,\,t)}=\max_{s\leq u\leq t}B^2_u$ and $0<s<t\leq T$.

\begin{prop}\label{B1_M2st}\textit{(Lee (2004))}
\begin{eqnarray}
\lefteqn{\mathbb{P}\{B^1_T\leq x,\,M^2_{(s,\,t)}\leq y\}}\nonumber\\
&=&\Phi_3\left(\frac{x-\mu_1 T}{\sigma_1\sqrt{T}},\,\frac{y-\mu_2 t}{\sigma_2\sqrt{t}},\,\frac{y-\mu_2 s}{\sigma_2\sqrt{s}};\,\rho\sqrt{\frac{t}{T}},\,\rho\sqrt{\frac{s}{T}},\,\sqrt{\frac{s}{t}}\right){}\label{M2st_a}\\
& &{}-e^{{2\,\mu_2\, y}/{\sigma_2^2}}\Phi_3\left(\frac{x-2\rho\frac{\sigma_1}{\sigma_2}y-\mu_1 T}{\sigma_1\sqrt{T}},\,\frac{-y-\mu_2 t}{\sigma_2\sqrt{t}},\,\frac{y+\mu_2 s}{\sigma_2\sqrt{s}};\,\rho\sqrt{\frac{t}{T}},\,-\rho\sqrt{\frac{s}{T}},\,-\sqrt{\frac{s}{t}}\right).\nonumber
\end{eqnarray}
\end{prop}  

\

\noindent \begin{dem} Let $\{Z_t\}_{t\in\mathbb{R}_+}$ be a stochastic process defined by
	$$Z_t=\frac{\sigma_2}{\sigma_1}B^1_t-\rho B^2_t\quad\forall t\in \mathbb{R}_+.$$
	It follows from the construction of $B^1$ and $B^2$ that the process $Z$ is a ${\cal{BM}}$ independent of $B^2$ with drift and variance parameters given by $\left(\frac{\sigma_2}{\sigma_1}\mu_1-\rho\mu_2\right)$ and $\sigma_2^2(1-\rho^2)$, respectively. Thus,
	\begin{eqnarray}
	\lefteqn{\mathbb{P}\{B^1_T\leq x,\,M^2_{(s,\,t)}\leq y\}}\nonumber\\
	&=&\mathbb{P}\{\frac{\sigma_1}{\sigma_2}\left(Z_T+\rho B^2_T\right)\leq x,\,M^2_{(s,\,t)}\leq y\}\nonumber\\
	&=&\int_{-\infty}^{+\infty}\mathbb{P}\{\rho B^2_T\leq \frac{\sigma_2}{\sigma_1}x-z,\,M^2_{(s,\,t)}\leq y|Z_T=z\}\mathbb{P}\{Z_T\in {\rm{d}}z\}\nonumber\\
	&=&\int_{-\infty}^{+\infty}\mathbb{P}\{\rho B^2_T\leq \frac{\sigma_2}{\sigma_1}x-z,\,M^2_{(s,\,t)}\leq y\}{}\label{M2st_0}\\
	& &{}\times\phi\left(\frac{z-\left(\frac{\sigma_2}{\sigma_1}\mu_1-\rho\mu_2\right)T}{\sigma_2\sqrt{(1-\rho^2)T}}\right)\frac{{\rm{d}}z}{\sigma_2\sqrt{(1-\rho^2)T}}.\nonumber
	\end{eqnarray}
	Define $z^*={2}x/{\sigma_1}-z$ and consider first the case where $\rho<0$.  The following probability has to be determined:
	\begin{eqnarray}
	\lefteqn{\mathbb{P}\{\rho B^2_T\leq \frac{\sigma_2}{\sigma_1}x-z,\,M^2_{(s,\,t)}\leq y\}}\nonumber\\
	&=&\mathbb{P}\{B^2_T\geq \frac{1}{\rho}z^*,\,M^2_{(s,\,t)}\leq y\}\nonumber\\
	&=&\mathbb{P}\{M^2_{(s,\,t)}\leq y\}-\mathbb{P}\{B^2_T\leq  \frac{1}{\rho}z^*,\,M^2_{(s,\,t)}\leq y\}\nonumber
	\end{eqnarray}
	\begin{eqnarray}
	&=&\left[\Phi_2\left(\frac{y-\mu_2 t}{\sigma_2\sqrt{t}},\,\frac{y-\mu_2 s}{\sigma_2\sqrt{s}};\,\sqrt{\frac{s}{t}}\right)\right.\nonumber\\
	& &{}-\left.\Phi_3\left(\frac{\frac{1}{\rho}z^*-\mu_2 T}{\sigma_2\sqrt{T}},\,\frac{y-\mu_2 t}{\sigma_2\sqrt{t}},\,\frac{y-\mu_2 s}{\sigma_2\sqrt{s}};\,\sqrt{\frac{t}{T}},\,\sqrt{\frac{s}{T}},\,\sqrt{\frac{s}{t}}\right)\right]{}\nonumber\\
	& &{}-e^{\frac{2\mu_2 y}{\sigma_2^2}}\left[\Phi_2\left(\frac{-y-\mu_2 t}{\sigma_2\sqrt{t}},\,\frac{y+\mu_2 s}{\sigma_2\sqrt{s}};\,-\sqrt{\frac{s}{t}}\right)\right.{}\nonumber\\
	& &{}-\left.\Phi_3\left(\frac{\frac{1}{\rho}z^*-2y-\mu_2 T}{\sigma_2\sqrt{T}},\,\frac{-y-\mu_2 t}{\sigma_2\sqrt{t}},\,\frac{y+\mu_2 s}{\sigma_2\sqrt{s}};\,\sqrt{\frac{t}{T}},\,-\sqrt{\frac{s}{T}},\,-\sqrt{\frac{s}{t}}\right)\right],\nonumber
	\end{eqnarray}
	which on applying equation (\ref{Phi_2_3}) becomes
	\begin{eqnarray}
	\lefteqn{\mathbb{P}\{\rho B^2_T\leq \frac{\sigma_2}{\sigma_1}x-z,\,M^2_{(s,\,t)}\leq y\}}\nonumber\\
	&=&\Phi_3\left(\frac{-(z^*-\rho\mu_2 T)}{\rho\sigma_2\sqrt{T}},\,\frac{y-\mu_2 t}{\sigma_2\sqrt{t}},\,\frac{y-\mu_2 s}{\sigma_2\sqrt{s}};\,-\sqrt{\frac{t}{T}},\,-\sqrt{\frac{s}{T}},\,\sqrt{\frac{s}{t}}\right){}\nonumber\\
	& &{}-e^{\frac{2\mu_2 y}{\sigma_2^2}}\Phi_3\left(\frac{-(z^*-\rho(2y+\mu_2 T))}{\rho\sigma_2\sqrt{T}},\,\frac{-y-\mu_2 t}{\sigma_2\sqrt{t}},\,\frac{y+\mu_2 s}{\sigma_2\sqrt{s}};\,-\sqrt{\frac{t}{T}},\,\sqrt{\frac{s}{T}},\,-\sqrt{\frac{s}{t}}\right)\nonumber
	\end{eqnarray}
	\begin{eqnarray}
	&=&\Phi_3\left(\frac{z^*-\rho\mu_2 T}{|\rho|\sigma_2\sqrt{T}},\,\frac{y-\mu_2 t}{\sigma_2\sqrt{t}},\,\frac{y-\mu_2 s}{\sigma_2\sqrt{s}};\,-\sqrt{\frac{t}{T}},\,-\sqrt{\frac{s}{T}},\,\sqrt{\frac{s}{t}}\right){}\label{M2st_2}\\
	& &{}-e^{\frac{2\mu_2 y}{\sigma_2^2}}\Phi_3\left(\frac{z^*-\rho(2y+\mu_2 T)}{|\rho|\sigma_2\sqrt{T}},\,\frac{-y-\mu_2 t}{\sigma_2\sqrt{t}},\,\frac{y+\mu_2 s}{\sigma_2\sqrt{s}};\,-\sqrt{\frac{t}{T}},\,\sqrt{\frac{s}{T}},\,-\sqrt{\frac{s}{t}}\right).\nonumber
	\end{eqnarray}
	
	Similary, when $\rho>0$, we have
	\begin{eqnarray}
	\lefteqn{\mathbb{P}\{\rho B^2_T\leq \frac{\sigma_2}{\sigma_1}x-z,\,M^2_{(s,\,t)}\leq y\}}\nonumber\\
	&=&\mathbb{P}\{B^2_T\leq \frac{1}{\rho}z^*,\,M^2_{(s,\,t)}\leq y\}\nonumber\\
	&=&\Phi_3\left(\frac{z^*-\rho\mu_2 T}{\rho\sigma_2\sqrt{T}},\,\frac{y-\mu_2 t}{\sigma_2\sqrt{t}},\,\frac{y-\mu_2 s}{\sigma_2\sqrt{s}};\,\sqrt{\frac{t}{T}},\,\sqrt{\frac{s}{T}},\,\sqrt{\frac{s}{t}}\right){}\label{M2st_1}\\
	& &{}-e^{\frac{2\mu_2 y}{\sigma_2^2}}\Phi_3\left(\frac{z^*-\rho(2y+\mu_2 T)}{\rho\sigma_2\sqrt{T}},\,\frac{-y-\mu_2 t}{\sigma_2\sqrt{t}},\,\frac{y+\mu_2 s}{\sigma_2\sqrt{s}};\,\sqrt{\frac{t}{T}},\,-\sqrt{\frac{s}{T}},\,-\sqrt{\frac{s}{t}}\right).\nonumber
	\end{eqnarray}
	
	On combining equations (\ref{M2st_2}) and (\ref{M2st_1}), we obtain the probability formula
	\begin{eqnarray}
	\lefteqn{\mathbb{P}\{\rho B^2_T\leq \frac{\sigma_2}{\sigma_1}x-z,\,M^2_{(s,\,t)}\leq y\}}\nonumber\\
	&=&\Phi_3\left(\frac{z^*-\rho\mu_2 T}{|\rho|\sigma_2\sqrt{T}},\,\frac{y-\mu_2 t}{\sigma_2\sqrt{t}},\,\frac{y-\mu_2 s}{\sigma_2\sqrt{s}};\,s(\rho)\sqrt{\frac{t}{T}},\,s(\rho)\sqrt{\frac{s}{T}},\,\sqrt{\frac{s}{t}}\right){}\label{M2st_3}\\
	& &{}-e^{\frac{2\mu_2 y}{\sigma_2^2}}\Phi_3\left(\frac{z^*-\rho(2y+\mu_2 T)}{|\rho|\sigma_2\sqrt{T}},\,\frac{-y-\mu_2 t}{\sigma_2\sqrt{t}},\,\frac{y+\mu_2 s}{\sigma_2\sqrt{s}};\,s(\rho)\sqrt{\frac{t}{T}},\,-s(\rho)\sqrt{\frac{s}{T}},\,-\sqrt{\frac{s}{t}}\right)\nonumber
	\end{eqnarray}
	where $s(\rho)=1$ if $\rho>0$ and $-1$ otherwise.
	
	In light of equation $(\ref{M2st_3})$, the result given in equation $(\ref{M2st_0})$ can be written as follows:
	\begin{eqnarray}
	\lefteqn{\mathbb{P}\{B^1_T\leq x,\,M^2_{(s,\,t)}\leq y\}}\nonumber\\
	&=&\int_{-\infty}^{+\infty}\Phi_3\left(\frac{z^*-\rho\mu_2 T}{|\rho|\sigma_2\sqrt{T}},\,\frac{y-\mu_2 t}{\sigma_2\sqrt{t}},\,\frac{y-\mu_2 s}{\sigma_2\sqrt{s}};\,s(\rho)\sqrt{\frac{t}{T}},\,s(\rho)\sqrt{\frac{s}{T}},\,\sqrt{\frac{s}{t}}\right){}\nonumber\\
	& &\times\phi\left(\frac{z-\left(\frac{\sigma_2}{\sigma_1}\mu_1-\rho\mu_2\right)T}{\sigma_2\sqrt{(1-\rho^2)T}}\right)\frac{{\rm{d}}z}{\sigma_2\sqrt{(1-\rho^2)T}}{}\nonumber\\
	& &{}-e^{\frac{2\mu_2 y}{\sigma_2^2}}\int_{-\infty}^{+\infty}\Phi_3\left(\frac{z^*-\rho(2y+\mu_2 T)}{|\rho|\sigma_2\sqrt{T}},\,\frac{-y-\mu_2 t}{\sigma_2\sqrt{t}},\,\frac{y+\mu_2 s}{\sigma_2\sqrt{s}};\,s(\rho)\sqrt{\frac{t}{T}},\,-s(\rho)\sqrt{\frac{s}{T}},\,-\sqrt{\frac{s}{t}}\right)\nonumber{}\\
	& &{}\times\phi\left(\frac{z-\left(\frac{\sigma_2}{\sigma_1}\mu_1-\rho\mu_2\right)T}{\sigma_2\sqrt{(1-\rho^2)T}}\right)\frac{{\rm{d}}z}{\sigma_2\sqrt{(1-\rho^2)T}}\nonumber\\
	&=&\Phi_3\left(\frac{x-\mu_1 T}{\sigma_1\sqrt{T}},\,\frac{y-\mu_2 t}{\sigma_2\sqrt{t}},\,\frac{y-\mu_2 s}{\sigma_2\sqrt{s}};\,\rho\sqrt{\frac{t}{T}},\,\rho\sqrt{\frac{s}{T}},\,\sqrt{\frac{s}{t}}\right){}\nonumber\\
	& &{}-e^{\frac{2\mu_2 y}{\sigma_2^2}}\Phi_3\left(\frac{x-2\rho\frac{\sigma_1}{\sigma_2}y-\mu_1 T}{\sigma_1\sqrt{T}},\,\frac{-y-\mu_2 t}{\sigma_2\sqrt{t}},\,\frac{y+\mu_2 s}{\sigma_2\sqrt{s}};\,\rho\sqrt{\frac{t}{T}},\,-\rho\sqrt{\frac{s}{T}},\,-\sqrt{\frac{s}{t}}\right)\nonumber
	\end{eqnarray}
where the last equality follows from Proposition~\ref{Int_Phi3} by letting  $a$ tend to $+\infty$.
\end{dem}

\begin{rem}When $\rho\rightarrow0$,
$$\mathbb{P}\{B^1_T\leq x,\,M^2_{(s,\,t)}\leq y\}\rightarrow \mathbb{P}\{B^1_T\leq x\}\mathbb{P}\{M^2_{(s,\,t)}\leq y\}.$$

Additionally, when $\mu_1=\mu_2$, $\sigma_1=\sigma_2$ and $\rho=1$, 
$$\mathbb{P}\{B^1_T\leq x,\,M^2_{(s,\,t)}\leq y\}=\mathbb{P}\{W^{(\mu_1,\,\sigma_1)}_T\leq x,\,M^{(\mu_1,\,\sigma_1)}_{(s,\,t)}\leq y\}.$$
\end{rem}

\section{A New Class of Bivariate Copulae}
In this section, several  bivariate copulae are constructed from the joint distribution functions  specified in the previous section.

In light of the invariance properties of copulae, we consider a ${\cal{BM}}$ with $\sigma=1$, since a $(\mu,\sigma)$-${\cal{BM}}$ can be derived from a $\left(\frac{\mu}{\sigma},\,1\right)$-${\cal{BM}}$ via a simple transformation (rescaling).

As well, it follows  from equations (\ref{FC_MBS}) and (\ref{F_Max}) that
$$F_{M_t}(a)=\mathbb{P}\{M_t\leq a\}=2\Phi\left(\frac{a}{\sqrt{t}}\right)-1,$$
and
\begin{eqnarray*}
F_{W_t,M_t}(x,a)&=&\mathbb{P}\{W_t\leq x,\,M_t\leq a\}\\
&=&\left\{\begin{array}{ll}
\Phi\left(\frac{x}{\sqrt{t}}\right)-\Phi\left(\frac{x-2a}{\sqrt{t}}\right) & \textrm{if }x\leq a\\
2\Phi\left(\frac{a}{\sqrt{t}}\right)-1 &\textrm{if } x>a.
\end{array}\right.
\end{eqnarray*}
Let 
$$F_{W_t}(x)=\mathbb{P}\{W_t\leq x\}=\Phi\left(\frac{x}{\sqrt{t}}\right)$$
be the marginal distribution of a standard ${\cal{BM}}$.

\begin{figure}[h!]
\begin{center}
\includegraphics[scale=0.45]{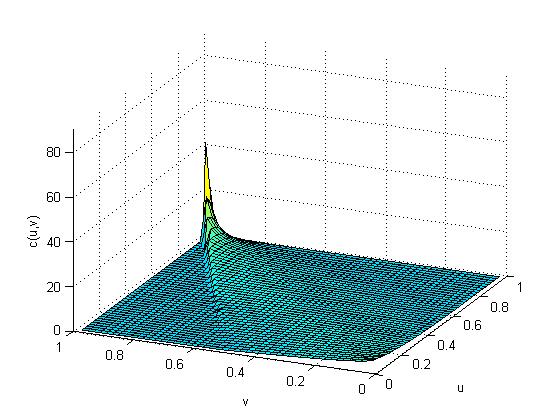}
\end{center}
\caption{Density of the copula generated by $W_t$ and $M_t$.}
\label{C_Mt}
\end{figure}

It follows from equation (\ref{C'=F}) that the copula $C_{W_t,\,M_t}(u,v)$ generated by a ${\cal{BM}}$ and its maximum is
\begin{equation}
C_{W_t,\,M_t}(u,v)=\left\{\begin{array}{ll}\label{C1}
u-\Phi\left(\Phi^{-1}(u)-2\Phi^{-1}\left(\frac{v+1}{2}\right)\right) & \textrm{if }u\leq \frac{v+1}{2}\\
v &\textrm{if } u>\frac{v+1}{2},
\end{array}\right.
\end{equation}
 its associated density function $c_{M_t}(u,\,v)$ being
\begin{eqnarray}
\lefteqn{c_{W_t,\,M_t}(u,\,v)=\frac{\partial^2}{\partial u\partial v}C_{W_t,\,M_t}(u,v)}\nonumber\\
&=&\frac{\left[2\Phi^{-1}\left(\frac{v+1}{2}\right)-\Phi^{-1}(u)\right]\phi\left(2\Phi^{-1}\left(\frac{v+1}{2}\right)-\Phi^{-1}(u)\right)}{\phi\left(\Phi^{-1}\left(\frac{v+1}{2}\right)\right)\phi\left(\Phi^{-1}(u)\right)}\label{c1}
\end{eqnarray}
whenever $u\leq \frac{v+1}{2}$, and zero otherwise.
This density is plotted in Figure~\ref{C_Mt}.
\\

Since  the copulae discussed in this paper involve  variables  that are not even interchangeable, they do not belong to the Archimedean class of copulae. Moreover, they clearly do not belong to the class of Gaussian copulae. They actually constitute a  new type of copulae whose distributions conglomerate in the neighborhood of the point (1,1) and, to a lesser extent, near the origin, the corresponding copula density functions being equal to zero beyond a certain treshold that is specified by a  relationship between the variables.

Let 
\begin{eqnarray*}
F_{W_t,M_t}(x,\,y;\,\mu)&=&\mathbb{P}\{W^{(\mu,\,1)}_t\leq x, M^{(\mu,\,1)}_t\leq y\}\\
&=&\left\{\begin{array}{ll}
\Phi\left(\frac{x-\mu t}{\sqrt{t}}\right)-e^{2\mu y}\Phi\left(\frac{x-2y-\mu t}{\sqrt{t}}\right) & \textrm{if $x\leq y$}\\
\Phi\left(\frac{y-\mu t}{\sqrt{t}}\right)-e^{2\mu y}\Phi\left(\frac{-y-\mu t}{\sqrt{t}}\right) & \textrm{if $x>y$},
\end{array}\right.
\end{eqnarray*}
and
\begin{eqnarray*}
F_{M_t}(y;\,\mu)&=&\mathbb{P}\{M_t^{(\mu,\,1)}\leq y\}\\
&=&\Phi\left(\frac{y-\mu t}{\sqrt{t}}\right)-e^{2\mu y}\Phi\left(\frac{-y-\mu t}{\sqrt{t}}\right),
\end{eqnarray*}
which are the distribution functions obtained in Proposition~\ref{FC_MAX} and equation (\ref{F_Max2}).

Let
$$F_{W_t}(x;\,\mu)=\mathbb{P}\{W_t^{(\mu,\,1)}\leq x\}=\Phi\left(\frac{x-\mu t}{\sqrt{t}}\right)$$
be the distribution function of a $(\mu,\,1)$-${\cal{BM}}$.
For $y>0$, the density function of $M_t^{(\mu,\,1)}$ is
\begin{eqnarray*}
f_{M_t}(y;\,\mu)&=&\frac{1}{\sqrt{t}}\phi\left(\frac{y-\mu t}{\sqrt{t}}\right)-e^{2\mu y}\left[2\mu\Phi\left(\frac{-y-\mu t}{\sqrt{t}}\right)-\frac{1}{\sqrt{t}}\phi\left(\frac{-y-\mu t}{\sqrt{t}}\right)\right].
\end{eqnarray*}

Therefore, the copula $C_{W_t,\,M_t}(u,v;\,\mu)$ generated by $W_t^{(\mu,\,1)}$ and $M_t^{(\mu,\,1)}$ is
\begin{eqnarray}
\lefteqn{C_{W_t,\,M_t}(u,v;\,\mu)}\nonumber\\
&=&\left\{\begin{array}{ll}
u-e^{2\mu \zeta(v)}\Phi\left(\Phi^{-1}(u)-\frac{2\zeta(v)}{\sqrt{t}}\right) & \textrm{if }u\leq \Phi\left(\frac{\zeta(v)-\mu t}{\sqrt{t}}\right)\label{C_Mt1}\\
v &\textrm{if } u>\Phi\left(\frac{\zeta(v)-\mu t}{\sqrt{t}}\right),
\end{array}\right.
\end{eqnarray}
and the corresponding density $c_{W_t,\,M_t}(u,v;\,\mu)$ is
\begin{eqnarray}
c_{W_t,\,M_t}(u,v;\,\mu)&=&\frac{2e^{2\mu \zeta(v)}\phi\left(\Phi^{-1}(u)-\frac{2\zeta(v)}{\sqrt{t}}\right)}{ f_{M_t}\left(\zeta(v);\,\mu\right)\phi\left(\Phi^{-1}(u)\right)}{}\label{cMt_pdf}\\
& &{}\times\left[\frac{1}{\sqrt{t}}\left(\frac{2\zeta(v)}{\sqrt{t}}-\Phi^{-1}(u)\right)-\mu\right]\quad\textrm{if }u\leq \Phi\left(\frac{\zeta(v)-\mu t}{\sigma\sqrt{t}}\right)\nonumber,
\end{eqnarray}
where $\zeta(v)=F_{M_t}^{-1}(v;\,\mu)$.

\begin{figure}[p]
\begin{center}
\includegraphics[scale=0.45]{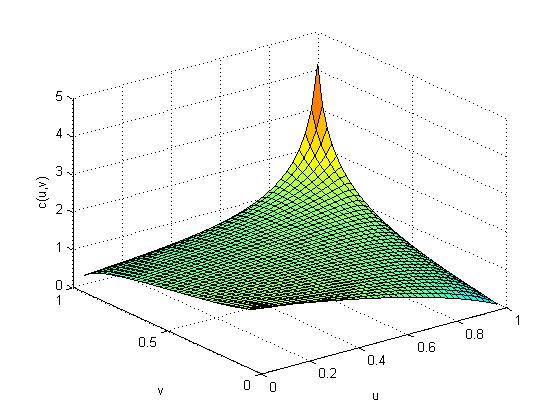}
\includegraphics[scale=0.45]{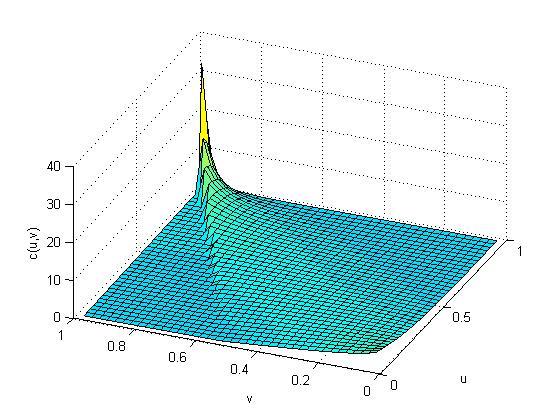}
\includegraphics[scale=0.45]{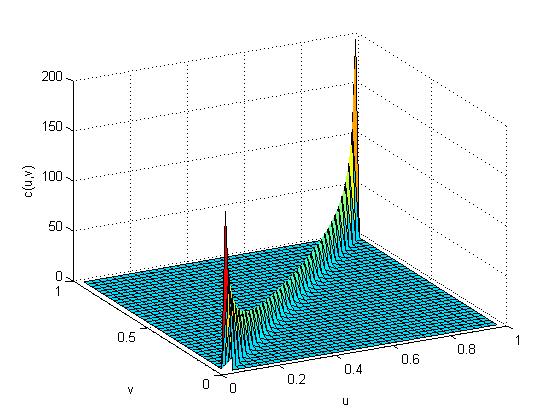}
\end{center}
\caption{Density functions of the copulae generated by $W_t^{(\mu,\,\sigma)}$ and $M_t^{(\mu,\,\sigma)}$ for increasing values of $\mu$.}
\label{C_Mt_mu}
\end{figure}

This density function appears in Figure~\ref{C_Mt_mu} for increasing values of $\mu$ ($\mu=-2$, $\mu=0$ and $\mu=10$, respectively). Clearly,  the strength of the dependence increases with $\mu$; additionally, as $\mu\rightarrow 0$, $C_{M_t}(u,\,v;\,\mu)\rightarrow C_{M_t}(u,\,v)$.

We know from Propositions~\ref{M_t_W_T} et \ref{M_st} that
\begin{eqnarray*}
F_{W_T,\,M_t}(x,\,y;\,\mu)&=&\mathbb{P}\{W_T^{(\mu,1)}\leq x,\,M_t^{(\mu,\,1)}\leq y\}\nonumber\\
&=&\Phi_2\left(\frac{x-\mu T}{\sqrt{T}},\,\frac{y-\mu t}{\sqrt{t}};\,\sqrt{\frac{t}{T}}\right){}\\
& &{}-e^{2\mu y}\Phi_2\left(\frac{x-2y-\mu T}{\sqrt{T}},\,\frac{-y-\mu t}{\sqrt{t}};\,\sqrt{\frac{t}{T}}\right),\nonumber
\end{eqnarray*}
\begin{eqnarray*}
F_{W_T,\,M_{(s,\,t)}}(x,\,y;\,\mu)&=&\mathbb{P}\{W_T^{(\mu,\,1)}\leq x,\,M_{(s,\,t)}^{(\mu,\,1)}\leq y\}\nonumber\\
&=&\Phi_3\left(\frac{x-\mu T}{\sqrt{T}},\,\frac{y-\mu t}{\sqrt{t}},\,\frac{y-\mu s}{\sqrt{s}};\,\sqrt{\frac{t}{T}},\,\sqrt{\frac{s}{T}},\,\sqrt{\frac{s}{t}}\right){}\\
& &{}-e^{2\mu y}\Phi_3\left(\frac{x-2y-\mu T}{\sqrt{T}},\,\frac{-y-\mu t}{\sqrt{t}},\,\frac{y+\mu s}{\sqrt{s}};\,\sqrt{\frac{t}{T}},\,-\sqrt{\frac{s}{T}},\,-\sqrt{\frac{s}{t}}\right),\nonumber
\end{eqnarray*}
and
\begin{eqnarray*}
F_{M_{(s,\,t)}}(y;\,\mu)&=&\mathbb{P}\{M_{(s,\,t)}^{(\mu,\,1)}\leq y\}\nonumber\\
&=&\Phi_2\left(\frac{y-\mu t}{\sigma\sqrt{t}},\,\frac{y-\mu s}{\sigma\sqrt{s}};\,\sqrt{\frac{s}{t}}\right){}\label{copule_1}\\
& &{}-e^{\frac{2\mu y}{\sigma^2}}\Phi_2\left(\frac{-y-\mu t}{\sigma\sqrt{t}},\,\frac{y+\mu s}{\sigma\sqrt{s}};\,-\sqrt{\frac{s}{t}}\right).\nonumber
\end{eqnarray*}

The copula $C_{W_T,\,M_t}(u,\,v;\mu)$ (resp. $C_{W_T,\,M_{(s,\,t)}}(u,\,v;\mu)$) describes the dependence structure induced by  $W_T^{(\mu,\,1)}$ and its maximum value on the time interval $[0,\,t]$(resp. $[s,\,t]$). Invoking (\ref{C'=F}), we obtain
\begin{eqnarray}
\lefteqn{C_{W_T,\,M_t}(u,\,v;\mu)}\nonumber\\
&=&\Phi_2\left(\Phi^{-1}(u),\,\frac{\zeta_1(v)-\mu t}{\sqrt{t}};\,\sqrt{\frac{t}{T}}\right){}\label{copule2}\\
& &{}-e^{2\mu \zeta_1(v)}\Phi_2\left(\Phi^{-1}(u)-\frac{2\zeta_1(v)}{\sqrt{T}},\,\frac{-\zeta_1(v)-\mu t}{\sqrt{t}};\,\sqrt{\frac{t}{T}}\right),\phantom{\frac{\zeta_2(v)+\mu s}{\sqrt{s}};\,-\sqrt{\frac{s}{T}},\,-\sqrt{\frac{s}{t}}}\nonumber
\end{eqnarray}
and
\begin{eqnarray}
\lefteqn{C_{W_T,\,M_{(s,\,t)}}(u,\,v;\mu)}\nonumber\\
&=&\Phi_3\left(\Phi^{-1}(u),\,\frac{\zeta_2(v)-\mu t}{\sqrt{t}},\,\frac{\zeta_2(v)-\mu s}{\sqrt{s}};\,\sqrt{\frac{t}{T}},\,\sqrt{\frac{s}{T}},\,\sqrt{\frac{s}{t}}\right){}\label{copule3}\\
& &{}-e^{2\mu \zeta_2(v)}\Phi_3\left(\Phi^{-1}(u)-\frac{2\zeta_2(v)}{\sqrt{T}},\,\frac{-\zeta_2(v)-\mu t}{\sqrt{t}},\,\frac{\zeta_2(v)+\mu s}{\sqrt{s}};\,\sqrt{\frac{t}{T}},\,-\sqrt{\frac{s}{T}},\,-\sqrt{\frac{s}{t}}\right),\nonumber
\end{eqnarray}
where $\zeta_1(v)=F^{-1}_{M_t}(v;\,\mu)$ and $\zeta_2(v)=F^{-1}_{M_{(s,\,t)}}(v;\,\mu)$.

Finally, consider $(B^1,\, B^2)$ a $\left(\mathbf{\mu^{*}},\mathbf{\Sigma}\right)$-${\cal{BM}}$ where
$$\mathbf{\mu^{*}}=(0,\,  \mu)'$$
and 
$$\mathbf{\Sigma}=\left(\begin{array}{cc}
1 & \rho\\
\rho & 1
\end{array}\right).$$
\begin{figure}[p]
	\begin{center}
		\includegraphics[scale=0.45]{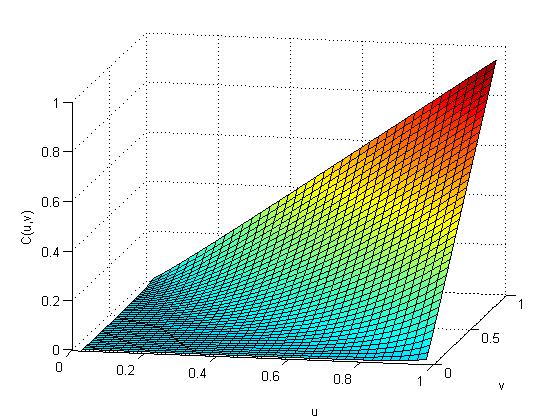}
		\includegraphics[scale=0.45]{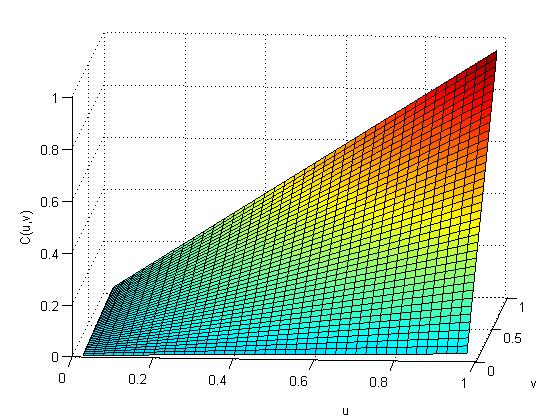}
		\includegraphics[scale=0.45]{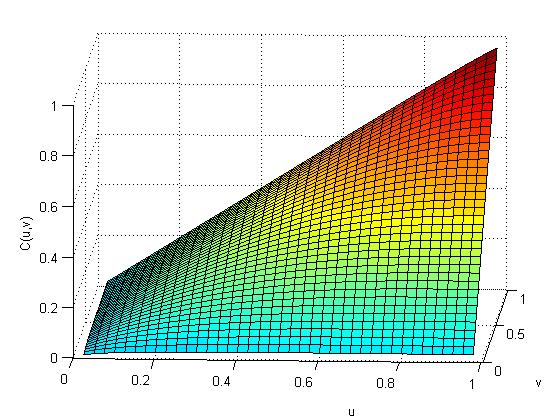}
	\end{center}
	\caption{Copulae generated by $B^1_T$ and $M^2_{(s,\,t)}$ for increasing values of $\rho$.}
	\label{C_M2st_B1T}
\end{figure}
The first ${\cal{BM}}$ has a zero drift because of the invariance property of copulae. Hence, from Proposition~\ref{B1_M2st}, we have
\begin{eqnarray*}
\lefteqn{F_{B^1_T,\,M^2_{(s,t)}}(x,y;\mu,\,\rho)}\nonumber\\
&=&\mathbb{P}\{B^1_T\leq x,\,M^2_{(s,\,t)}\leq y\}\nonumber\\
&=&\Phi_3\left(\frac{x}{\sqrt{T}},\,\frac{y-\mu t}{\sqrt{t}},\,\frac{y-\mu s}{\sqrt{s}};\,\rho\sqrt{\frac{t}{T}},\,\rho\sqrt{\frac{s}{T}},\,\sqrt{\frac{s}{t}}\right){}\\
& &{}-e^{2\mu y}\Phi_3\left(\frac{x-2\rho y}{\sqrt{T}},\,\frac{-y-\mu t}{\sqrt{t}},\,\frac{y+\mu s}{\sqrt{s}};\,\rho\sqrt{\frac{t}{T}},\,-\rho\sqrt{\frac{s}{T}},\,-\sqrt{\frac{s}{t}}\right).\nonumber
\end{eqnarray*}
Let us now denote by
$$F_{B^1_T}(x)=\Phi\left(\frac{x}{\sqrt{T}}\right)$$
the distribution function of the first ${\cal{BM}}$ and by $F_{M_{(s,\,t)}^2}(y;\,\mu)$ the distribution function of $M_{(s,\,t)}^2$ where $F_{M_{(s,\,t)}^2}(y;\,\mu)=F_{M_{(s,\,t)}}(y;\,\mu)$ for all $y>0$.

The bivariate copula $C_{B^1_T,\,M^2_{(s,\,t)}}(u,v;\,\mu,\,\rho)$ generated by $B^1_T$ and $M^2_{(s,\,t)}$ is then defined by
\begin{eqnarray}
\lefteqn{C_{B^1_T,\,M^2_{(s,\,t)}}(u,v;\,\mu,\,\rho)}\nonumber\\
&=&\Phi_3\left(\Phi^{-1}(u),\,\frac{\zeta(v)-\mu t}{\sqrt{t}},\,\frac{\zeta(v)-\mu s}{\sqrt{s}};\,\rho\sqrt{\frac{t}{T}},\,\rho\sqrt{\frac{s}{T}},\,\sqrt{\frac{s}{t}}\right){}\label{Marie}\\
& &{}-e^{2\mu \zeta(v)}\Phi_3\left(\Phi^{-1}(u)-\frac{2\rho \zeta(v)}{\sqrt{T}},\,\frac{-\zeta(v)-\mu t}{\sqrt{t}},\,\frac{\zeta(v)+\mu s}{\sqrt{s}};\,\rho\sqrt{\frac{t}{T}},\,-\rho\sqrt{\frac{s}{T}},\,-\sqrt{\frac{s}{t}}\right)\nonumber
\end{eqnarray}
where $\zeta(v)=F^{-1}_{M^2_{(s,\,t)}}(v;\,\mu)$.
This copula contains all the copulae considered in this section. Indeed, when $\rho$ ten{\rm{d}}s to $1$, the copula specified by equation (\ref{Marie}) converges to that generated by a ${\cal{BM}}$ with drift $\mu$ at $T$ and its own maximum on the interval $[s,\,t]$, which is given in equation (\ref{copule3}). From this result, we obtain the copula given in equation (\ref{copule2}) by letting $s$ tend to $0$. Finally, as $\rho=1$, $s\rightarrow 0$ and $t\rightarrow T$, the copula specified by equation (\ref{Marie}) converges to that given in (\ref{C_Mt1}).

The copula  $C_{B^1_T,\,M^2_{(s,\,t)}}(u,\,v;\,\mu)$ is plotted in Figure~\ref{C_M2st_B1T} for $\rho=-0.99$, $\rho=0$ and $\rho=0.99$.

Note that when $\rho=0$, $C_{B^1_T,\,M^2_{(s,\,t)}}(u,\,v;\mu,\,\rho)=C_I(u,\,v)$ where $C_I$ is the independent copula defined by $C_I(u,\,v)=uv$ for all $(u,\,v)\in \mathbf{I}^2$.\\

The proposed copulae are  applicable to certain bivariate data sets for which one of the variables involves maxima.  Such observations occur for instance in hydrology,  meteorology and financial modeling.

\

\ 

\noindent {\bf Acknowledgements}\\

The financial support of the Natural Sciences and Engineering Research Council of Canada is gratefully acknowledged by the first and third authors. Thanks are also due to Arthur Charpentier, Jean-Fran{\c{c}}ois Plante and Bruno R{\'e}millard for their comments on an initial draft of the paper.

\

\noindent{\bf Orcid}\\

Serge B Provost ID  0000-0002-2024-0103

\

\end{document}